\documentclass[12pt]{amsart}
\usepackage{amssymb}
\usepackage{verbatim}
\usepackage{amsmath}

\theoremstyle{plain}
\newtheorem{theorem}{Theorem}[section]
\newtheorem{corollary}[theorem]{Corollary}
\newtheorem{proposition}[theorem]{Proposition}
\newtheorem{lemma}[theorem]{Lemma}
\newtheorem{conclusion}{Conclusion}

\newtheorem{fact}{Fact}

\theoremstyle{definition}

\theoremstyle{remark}
\newtheorem*{remark}{Remark}
\newtheorem*{notation}{Notation}

\newcounter{numb}
\newcounter{numb2}
\newenvironment{theoremlist}
   {\begin{list}
        {\arabic{numb2})}
    {\setlength{\topsep}{2 pt}
     \setlength{\itemsep}{2 pt} 
     \setlength{\leftmargin}{20 pt}
     \setlength{\rightmargin}{0 pt}
     \setlength{\itemindent}{0 pt}
      \usecounter{numb2} }
     }
   {\end{list}}

\newcommand{\ds}{\displaystyle}
\newcommand{\cstar}{\ensuremath{\text{C}^{*}}\nobreakdash-\hspace{0 pt}}
\renewcommand{\star}{\ensuremath{{}^{*}}\nobreakdash-\hspace{0 pt}}
\newcommand{\BB}{\mathcal{B}}
\newcommand{\DD}{\mathcal{D}}
\newcommand{\HH}{\mathcal{H}}
\newcommand{\II}{\mathcal{I}}
\newcommand{\JJ}{\mathcal{J}}
\newcommand{\KK}{\mathcal{K}}

\newcommand{\MM}{\mathcal{M}}
\newcommand{\NN}{\mathcal{N}}
\newcommand{\alg}{\operatorname{Alg}}
\newcommand{\lat}{\operatorname{Lat}}
\newcommand{\supp}{\operatorname{supp}}
\newcommand{\suc}{\operatorname{succ}}
\newcommand{\pred}{\operatorname{pred}}

\newcommand{\Ran}{\operatorname{Ran}}

\newcommand{\m}{\operatorname{m}}
\newcommand{\sptr}{\ensuremath{(X,P,G)}}
\newcommand{\s}{\sigma}
\newcommand{\orb}{\text{orb}}
\newcommand{\bara}{\bar{\alpha}}
\newcommand{\barb}{\bar{\beta}}
\newcommand{\ehat}{\hat{e}}
\newcommand{\fhat}{\hat{f}}
\newcommand{\uhat}{\hat{u}}
\newcommand{\vhat}{\hat{v}}
\newcommand{\what}{\hat{w}}

\newcommand{\sab}{\s_{a,b}}
\newcommand{\tab}{\tau_{a,b}}
\newcommand{\nsa}{non-self-adjoint}
\newcommand{\Nsa}{Non-self-adjoint}
\newcommand{\np}{n\nobreakdash-\hspace{0 pt}primitive}
\newcommand{\dirlim}{\varinjlim}

\newcommand{\wh}{\widehat}

\begin{document}
\addtolength{\textwidth}{25 pt}
\title[Nest Representations]{ Nest Representations of TAF Algebras }
\author{Alan Hopenwasser}
\address{Department of Mathematics\\
        University of Alabama\\
        Tuscaloosa, AL 35487, USA}
\email{ahopenwa@euler.math.ua.edu}

  \author{Justin~R. Peters}
  \address{Department of Mathematics\\
        Iowa State University\\
        Ames, IA, USA}
  \email{peters@iastate.edu} 

  \author{Stephen~C. Power}
  \address{Department of Mathematics and Statistics\\
        Lancaster University\\
        Lancaster, LA1 4YF, UK}
  \email{s.power@lancaster.ac.uk}

\thanks{A.~Hopenwasser was supported in part by an EPSRC visiting fellowship
 to Lancaster University.}
\thanks{Justin.~R.~Peters was partially supported by a grant from the
  National Science Foundation.}
\thanks{EPSCoR of Alabama provided a travel grant to S.~C.~Power for
a visit the University of Alabama.}

 \keywords{Nest representation, meet irreducible ideal, strongly
maximal TAF algebra.}
 \subjclass{47D25}
 \date{August 28, 1999}

\begin{abstract}
A nest representation of a strongly maximal TAF algebra $A$
with diagonal $D$ 
is a representation $\pi$ for which $\lat \pi(A)$ is totally
ordered.  We prove that $\ker \pi$ is a meet irreducible ideal
 if the spectrum of $A$ is totally 
ordered or if (after an appropriate similarity)
the von Neumann algebra $\pi(D)''$
 contains an atom.
 \end{abstract}
\maketitle

\section{Introduction}\label{s:intro}

Irreducible \star representations and their kernels,
the primitive ideals, play a fundamental role in the theory
of \cstar algebras.  For example, the structure of the lattice
of all ideals in a \cstar algebra is determined by the space of
primitive ideals with the hull-kernel topology (there is
a bijection between the lattice of ideals and the lattice of
closed subsets of the primitive ideal space) and every ideal
is the intersection of those primitive ideals which contain it.

Recent work by Michael Lamoureux \cite{mpl93, mpl96, mpl97} has
shown that a similar situation prevails in a number of
\nsa\  operator algebra settings.  One motivation for
investigating \nsa\  algebras arises from dynamical systems.
While \cstar crossed products constructed from dynamical
systems are very useful in the study of dynamical systems,
some essential information may be lost.  (There are different
dynamical systems which give rise to isomorphic \cstar crossed
products.)  In the case of free, discrete systems, the remedy is
to look instead at the semi-crossed product (a \nsa\  algebra)
associated with the system (see \cite{wakj69,jp84, scp92a}). 
 In this context there
is a bijection between the isomorphism classes of (free, discrete)
dynamical systems and their associated semi-crossed products.  
This correspondance is established via an analysis of the space
of ideals of the operator algebra.

In attempting to extend this program to other dynamical systems,
Lamoureux drew attention to nest representations and their
kernels (which he called \textit{\np\  ideals}).  A
\textit{nest representation} $\pi$ of an operator
algebra $A$ is simply a continuous algebra homomorphism
of $A$ into some $\BB(\HH)$ with the property that the lattice
of projections invariant under $\pi(A)$ is totally ordered.
Nest representations do not arise naturally as a concept in the
\cstar algebra context.  For one thing, most representations arising
in \cstar algebra theory are \star representations.  
A projection is invariant under a \star representation
 if, and only if, it is reducing; 
consequently, for \star representations,
the family of nest representations reduces to the family of 
irreducible representations.  Even if one looks at representations
which are not \star representations, the \np\  ideals in a 
\cstar algebra are just the primitive ideals.  (Nest representations
are necessarily topologically cyclic -- this is valid for
nest representations of any Banach algebra with a bounded approximate
identity \cite{mpl93} -- and Haagerup \cite{uh83} has shown that
a cyclic representation of a \cstar algebra is similar to
a \star representation.)

In the \nsa\  context, the story is quite different.  \Nsa\  algebras
may lack primitive ideals, but \np\  ideals generally abound.  Lamoureux
has shown that, at least in the presence of certain hypotheses,
\np\  ideals in semi-crossed products provide information about
arc spaces in dynamical systems, so that the semi-crossed product
essentially determines the dynamical system.  He
 has also shown that, in a variety of contexts, the \np\  ideals
carry a topology and that the lattice of closed sets in this topology
is isomorphic to the lattice of ideals in the original algebra.
Also, every ideal is equal to the intersection of the \np\  ideals
which contain it.  Thus, for general operator algebras, \np\  ideals
are often a suitable replacement for the primitive ideals of 
\cstar algebra theory.  

The topology which Lamoureux puts on the \np\  ideals is, of course,
the hull-kernel topology.  In order to show that the hull-kernel 
operation yields a topology, Lamoureux uses a technical property for
ideals which is related to meet irreducibility and which implies
meet irreducibility.  (An ideal $\II$ is 
\textit{meet irreducible} if, for
any ideals $\JJ$ and $\KK$, $\II = \JJ \cap 
\KK \Longrightarrow \II = \JJ \text{ or } \II = \KK$.)
  Meet irreducibility arises in \cite{mpl96} in connection
with semi-crossed products and dynamical systems.

The following four results from \cite{mpl97} indicate that
meet irreducibility will be closely related to 
n-primitivity in diverse operator algebra contexts. 
\begin{list}{(\arabic{numb})}
    {\setlength{\topsep}{12 pt}
     \setlength{\itemsep}{10 pt} 
     \setlength{\leftmargin}{24 pt}
     \setlength{\rightmargin}{0 pt}
     \setlength{\itemindent}{0 pt}
      \usecounter{numb} }
\item Let $\II$ be a closed, two-sided ideal in a separable
\cstar algebra.  Then $\II$~is \np\  $\Longleftrightarrow$
 $\II$~is primitive $\Longleftrightarrow$ 
$\II$~is prime $\Longleftrightarrow$
$\II$~is meet irreducible.  (Some of this has, of course, been
known for a long while.)
\item Let $\II$ be a closed, two-sided ideal in
$\alg(\NN)\cap{\KK}$, where $\NN$ is a nest of closed subspaces
in some Hilbert space $\HH$, $\alg(\NN)$ is the associated
nest algebra consisting of all operators which leave invariant
each subspace in $\NN$, and $\KK$ is the algebra of all compact
operators acting on $\HH$.  Then, $\II$~is \np\  $\Longleftrightarrow$
$\II$~is meet irreducible $\Longleftrightarrow$ $\II$~is the kernel of
the compression map of $\alg(\NN) \cap \KK$ to some interval
from $\NN$.
\item Let $\II$ be a closed ideal in the disk algebra, $A(\mathbb{D})$,
the algebra of continuous functions on the unit disk of $\mathbb{C}$
which are analytic in the interior.  Then, $\II$~is \np\ 
$\Longleftrightarrow$ $\II$~is meet irreducible $\Longleftrightarrow$
$\II$~is either primary or zero.  (In this context, it is not
true that every ideal is the intersection of meet irreducible
ideals.)
\item Let $A=T_{n_1} \oplus \dots \oplus T_{n_k}$ be a direct
sum of upper triangular $n_j$ by $n_j$ matrix algebras and let
$\II$ be a two-sided ideal in $A$.  Then, $\II$~is \np\  
$\Longleftrightarrow$ $\II$~is meet irreducible.
\end{list}

Meet irreducible ideals in the context of triangular AF algebras
(TAF algebras) were investigated in \cite{dhhls}.  In particular,
it was proven that every meet irreducible ideal in a strongly
maximal TAF algebra is \np\ \cite[Theorem 2.4]{dhhls}.  Once again,
the lattice of all ideals is isomorphic to the lattice of closed sets
in the space of meet irreducible ideals with the hull-kernel topology
and every ideal is an intersection of meet irreducible ideals.  
The converse, ``every \np\ ideal is meet irreducible,'' was left
open, however; it is the purpose of this note to investigate
this converse in the context of strongly maximal TAF algebras.
In section~\ref{s:totor} we describe a broad class of algebras
(those characterized by ``totally ordered spectrum'') for which
the converse always holds.  In section~\ref{s:atomic} we show that
the converse is valid for any nest representation which is similar
to a representation $\pi$ which is a \star representation on the
diagonal $D$ of the algebra and for which the von Neumann algebra
generated by $\pi(D)$  contains an atom. 
 These two sections may be read independently
of each other.

\subsection{Notation} \label{ss:not}
We now establish notation and terminology.  Let $B$ be an
AF \cstar algebra and let $D$ be a canonical masa in $B$.
This implies that there is a sequence of finite dimensional
\cstar algebras $B_i$ and embeddings $\phi_i \colon B_i 
\rightarrow B_{i+1}$
 such that $\ds B = \dirlim (B_i, \phi_i)$ and 
$\ds D = \dirlim (D_i, \phi_i)$, where each
$D_i = D \cap B_i$ is a masa in $B_i$ and each
$\phi_i$ maps the normalizer of $D_i$ into the
normalizer of $D_{i+1}$.
In particular, the $D$-normalizing partial isometries in $B$
have a linear span which is dense in $B$
($w$ \textit{normalizes}
$D$ if $wDw^* \subseteq D$ and $w^*Dw \subseteq D$;
the set of $D$-normalizing partial isometries in $B$
is denoted by $N_D(B)$).
A TAF \textit{algebra} $A$ with diagonal $D$ is a subalgebra
of $B$  such that $A \cap A^* = D$.  It follows
that $\ds A = \dirlim (A_i, \phi_i)$, where 
$A_i = A \cap B_i$, for all $i$. Each $A_i$ is
necessarily triangular in $B_i$ with diagonal $D_i$;
if, in addition, each $A_i$ is maximal triangular
then we say that $A$ is \textit{strongly maximal}.
This is equivalent to requiring that $A+A^*$ be
dense in $B$.  

AF \cstar algebras are groupoid \cstar algebras
with groupoids which are especially tractable:
the groupoids are topological equivalence relations.
As such, the \cstar algebra can be identified with
an algebra of functions on the groupoid.  Subalgebras
such as TAF algebras and ideals determine and are
determined by appropriate substructures of the groupoid. 
This provides a coordinatization for TAF algebras and
their ideals.

The following is a brief sketch of this coordinatization
for strongly maximal TAF algebras.  For a more thorough
description, see \cite{scp92bk} or \cite{ms89}.  Let
$(D, A, B)$ be a triple in which $B$ is an AF
\cstar algebra, $D$ is a canonical masa, and
$A$ is a strongly maximal TAF subalgebra of $B$
with diagonal $D$.  We need to describe the
\textit{spectral triple}, $(X, P, G)$ associated
with $(D,A,B)$.  The first ingredient $X$ is the
usual spectrum of the abelian \cstar algebra $D$
(so that $D \cong C(X)$).  Note that, in the present
context, $X$ will be a Cantor space. 
 The algebra $B$ is generated
by partial isometries which normalize $D$. 
Each normalizing partial isometry induces a partial 
homeomorphism of $X$ into itself (a homeomorphism of
a closed subset of $X$ onto another closed subset).  
The union of the graphs of these homeomorphisms is
an equivalence relation;  when this set is topologized so that
the graph of each normalizing partial isometry is open and
closed, one obtains the groupoid $G$.  The spectrum $P$ of
$A$ is the union of the graphs of the normalizing 
partial isometries which lie in $A$.

If $x \in X$, the equivalence class in $G$ which
contains $x$ is referred to as the \textit{orbit}  of
$x$ (since it consists of all the images of $x$ under
homeomorphisms associated with $D$-normalizing partial 
isometries) and is denoted by $\orb_x$.
When $(x,y) \in P$, we shall often write $x \preceq y$;
when $A$ is strongly maximal, $P$ induces a total order on
each orbit.

Now suppose that $\pi$ is a nest representation of
a TAF algebra $A$.    Since $D$ is an abelian \cstar algebra,
results in \cite{rvk55} imply that the restriction of
$\pi$ to $D$ is similar to a \star representation.  (For the 
limited domain in which we need this theorem, abelian 
\cstar algebras, Kadison attributes this fact to unpublished
1952 lecture notes of Mackey.) Thus any nest representation
of a TAF algebra is similar to another nest representation
whose restriction to the diagonal is a \star representation.
Accordingly, we assume throughout this paper that any nest
representation $\pi$ of a TAF algebra $A$ acts as a
\star representation on the diagonal, $D$.

\section{Totally Ordered Spectrum}\label{s:totor}

In this section, $A$ will be a strongly maximal
 TAF algebra whose \cstar envelope $B$
is simple.  This implies that the orbit of each element
in $X$ is dense in $X$.
  The diagonal of $A$ is denoted by $D$.
The spectral triple for $(D,A,B)$ is $\sptr$.
We assume that $X$ has a total order $\preceq$
which agrees on each equivalence class with the total
order induced by $P$.
$X$ has a minimal and a maximal element, which we
denote by $0$ and $1$, respectively.
If $a \in X$, we say that $a$ has a \textit{gap above}
if $a$ has an immediate successor and that $a$ has 
a \textit{gap below} if $a$ has an immediate predecessor.

Elements of $B$ will be viewed as continuous functions
in $C_0(G)$; elements of $A$ are continuous functions whose
supports lie in $P$.  Of course, not all elements of $C_0(G)$ 
are elements of $B$, but all those with compact support
certainly are in $B$.  When elements of $B$ are viewed
as functions, the multiplication is given by a convolution
formula with respect to a Haar system consisting of
counting measure on each equivalence class from $G$.

If $Y$ is a clopen subset of $X$, then $E_Y$ denotes the
projection in $D$ corresponding to $Y$ (i.e., the characteristic
function of $Y$).

 In the situation where $X$ carries a total order compatible with $P$,
the meet irreducible ideal sets have been completely described in
Theorem 3.1 in \cite{dhhls}. For the convenience of the reader, we
restate that theorem here.  For each pair of elements $a$ and $b$
in $X$, define two subsets of $P$:
\begin{align*}
\sab &= \{ (x,y) \in P \mid x \prec a \text{ or }
  b \prec y \}, \\
\tab &= \sab \cup \{ (a,b) \}.
\end{align*}
The set $\sab$ is an ideal set in $P$; the set $\tab$ 
is an ideal set provided that $(a,b) \in P$ and $\tab$
is an open subset of $P$.  (We always assume that
these two conditions hold.)

\begin{theorem}\label{t:misets}
Assume that $A$ is a strongly maximal TAF algebra with
simple \cstar envelope and totally ordered spectrum.
The following is a complete list of all the meet irreducible
ideal sets in $P$:
\begin{theoremlist}
\item $\sab$, if $(a,b) \in P$;
\item $\sab$, if $(a,b) \notin P$ and either $a$ has
no gap above or $b$ has no gap below;
\item $\tab$, where either $a$ has no gap above or $b$
has no gap below.
\end{theoremlist}
\end{theorem}

\begin{theorem}\label{t:miker}
Let $\pi$ be a continuous nest representation of $A$, 
where $A$ is a strongly maximal TAF algebra with 
simple \cstar envelope and totally
ordered spectrum.  
Then the kernel of $\pi$ is a meet irreducible ideal
in $A$.
\end{theorem}

 \begin{proof}
Let $\pi$ be a continuous nest representation of $A$
acting on the Hilbert space $\HH$.  Let $\s$ be the
ideal set in $P$ which corresponds to the ideal
$\ker \pi$.  
Through a series of facts, we will show that $\s$
is one of the meet irreducible ideal sets listed
in Theorem~\ref{t:misets}; consequently,  $\ker \pi$
is meet irreducible.

\begin{fact}\label{f:fact1} Suppose that $a$ has a gap above or,
equivalently, that $E_{[0,a]}$ is a projection in $A$ not
equal to the identity.
Then $E_{[0,a]}$ is
invariant for $A$ and hence $\pi(E_{[0,a]}) \in \lat{\pi(A)}$.
\end{fact}

\begin{proof}
Let $f \in A$; view $f$ as a $C_0$ function on $G$.  For
any $(x,y) \in G$,

 \begin{align*}
 f E_{[0,a]} (x,y) &= \sum_{z \in \orb_x} f(x,z) E_{[0,a]} (z,y) \\
  &= \begin{cases}
        f(x,y), & \text{if } y \preceq a \\
        0,  & \text{otherwise}
     \end{cases} \\
  &= E_{[0,a]}(x,x) f(x,y) E_{[0,a]} (y,y) \\ 
  &= E_{[0,a]} f E_{[0,a]} (x,y).
 \end{align*}
We use the fact that $(x,y) \in \supp f$ implies that $x \preceq y$ in $X$.
\end{proof}

\begin{fact}\label{f:fact2}
If $a$ has a gap above and if there is a point 
$(a,c) \in P \setminus \s$, then $\pi(E_{[0,a]}) \ne 0$
 and $\pi(E_{[0,a]})$ is a non-trivial invariant
projection for $\pi(A)$.  
\end{fact}

\begin{remark} If we were not assuming that $\pi$ is
normalized to a \star representation of $D$, then
$\pi(E_{[0,a]})$ would be  a non-trivial
idempotent whose range is an invariant subspace for 
$\pi(A)$.
\end{remark}

\begin{proof}
Since $(a,a) \circ (a,c) = (a,c)$ and $\s$ is an ideal
set, $(a,a) \in P \setminus \s$.  If $\pi(E_{[0,a]}) = 0$,
then $E_{[0,a]} \in \ker{\pi}$; hence, $\supp {E_{[0,a]}}
\subseteq \s$.  Thus $(a,a) \in \s$, a contradiction. 
\end{proof}

\begin{fact}\label{f:fact3}
If $b$ has a gap below and if there is a point
$(c,b) \in P \setminus \s$, then $\pi(E_{[b,1]}) \ne 0$.
\end{fact}

\begin{proof}
Essentially the same as for Fact~\ref{f:fact2}:
$(c,b) \circ (b,b) = (c,b)$, so $(b,b) \in P \setminus \s$.
\end{proof}

\begin{fact}\label{f:fact4}
Assume $a$ has a gap above, $b$ has a gap below, 
$\pi(E_{[0,a]}) \ne 0$, and $\pi(E_{[b,1]}) \ne 0$.
If $E_{[0,a]} f E_{[b,1]} \in \ker{\pi}$ for all
 $f \in N_D(A)$ (or, for all $f$ in a set with 
linear span dense in $A$), then $\pi$ is not a 
nest representation.
\end{fact}

\begin{proof}
Since $\ker{\pi}$ is closed, $E_{[0, a]} f E_{[b, 1]}
 \in \ker{\pi}$, for all $f \in A$. As a consequence, we
have $\pi(E_{[0, a]}) \pi(f) \pi(E_{[b, 1]}) = 0 $,
for all $f \in A$.  Let $\KK_1$ be the range
of $\pi(E_{[0, a]})$ and $\KK_2$ be the norm closure
of $\pi(A) \pi(E_{[b, 1]}) \HH$.  Then $\KK_1$ and
$\KK_2$ are non-zero invariant subspaces for $\pi(A)$
and $\KK_1 \cap \KK_2 = (0)$.  So $\lat \pi(A)$ is not
totally ordered by inclusion; $\pi$ is not a nest 
representation.    
\end{proof} 

In what follows we use the standard identification of
$X$ with the diagonal of $P$,
i.e. with $\{ (x,x) \mid x \in X \}$.

\begin{fact}\label{f:fact5}
If $\pi$ is a nest representation, then
$(P \setminus \s) \cap X$ is an interval in $X$.
\end{fact}

\begin{proof}
Suppose not.  Then there are three points $a,b,c$
with $a \prec b \prec c $ in $X$
such that $(a,a) \in P \setminus \s$, $(b,b) \in \s$, and
$(c,c) \in P \setminus \s$.  It follows that if $(x,y) \in P$
with $x \preceq b$ and $b \preceq y$, then $(x,y) \in \s$.
(Since $\s$ is open there is a neighborhood of $(b,b)$ which
is contained in $\s$; use the assumption that all orbits are
dense and the ideal set property for $\s$.)

Choose $\alpha$ and $\beta$ so that $a \preceq \alpha \prec b$,
$b \prec \beta \preceq c$, $\alpha$ has a gap above, and $\beta$
has a gap below.  Then $(a,a) \in \supp E_{[0,\alpha]}$;
hence $\pi(E_{[o, \alpha]}) \ne 0$ and the range of
$\pi(E_{[0, \alpha]})$ is in $\lat \pi(A)$.
Also, $(b,b) \in \supp E_{[\beta, 1]}$, so 
$\pi(E_{[\beta, 1]}) \ne 0$.  For any $g \in A$,
\begin{align*}
(x,y) \in \supp E_{[0, \alpha]} g E_{[\beta, 1]}
  &\Longrightarrow  x \preceq \alpha \prec b \text{ and }
   b \prec \beta \preceq y \\
  &\Longrightarrow (x,y) \in \s.
\end{align*}
Consequently,  $E_{[0, \alpha]} g E_{[\beta, 1]} \in \ker \pi$.
Fact~\ref{f:fact4} now implies that $\pi$ is not a nest
representation, contradicting the hypothesis.
\end{proof}

Assume $\pi$ is a nest representation and that
$a$ is the left endpoint of $(P \setminus \s) \cap X$
and $b$ is the right endpoint of $(P \setminus \s)
 \cap X$.  Each of $a$ and $b$ may or may not be 
elements of $(P \setminus \s) \cap X$.  However,
if $a$ has a gap above, then without loss of
generality, we  may assume that
$a \in (P \setminus \s) \cap X$ (simply replace $a$
by its immediate successor, if necessary).  Similarly,
if $b$ has a gap below, we may assume that
$b \in (P \setminus \s) \cap X$.

\begin{fact}\label{f:fact6}
If $a \prec \alpha \preceq \beta \prec b$ and 
$(\alpha, \beta) \in P$, then $(\alpha, \beta)
 \in P \setminus \s$.
\end{fact}

\begin{proof}
Suppose that $a \prec \alpha \preceq \beta \prec b$,
$(\alpha, \beta) \in P$, and $(\alpha, \beta) \in \s$.
Choose $\bara$ and $\barb$ so that 
$a \prec \bara \preceq \alpha \preceq \beta \preceq \barb 
\prec b$, $\bara$ has a gap above, and $\barb$ has a gap
below. [Exceptions: if $\alpha = \suc a$, choose
$\bara = a$ and note that this element is not in $\s \cap X$;
if $\beta = \pred b$, choose $\barb = b$ and note that this
is not in $\s \cap X$.]

It follows that $\pi(E_{[0, \bara]})$ is non-zero and has
range in $\lat \pi(A)$ and that $\pi(E_{[\barb, 1]}) \ne 0$.
Let $g \in A$ be arbitrary.  If $(x,y)$ is in
$\supp E_{[0, \bara]} g E_{[\barb, 1]}$, then 
$x \preceq \bara \preceq \alpha$ and $\beta \preceq \barb
\preceq y$; hence $(x,y) \in \s$.  Thus 
$E_{[0, \bara]} g E_{[\barb, 1]} \in \ker \pi$, for all $g$.
Fact~\ref{f:fact4} now implies that $\pi$ is not a nest
representation, a contradiction.
\end{proof}

\begin{conclusion}\label{c:concA}
If $\pi$ is a nest representation and $a$ and $b$
are the endpoints of the interval $(P \setminus \s) \cap X$, then
\[
\{ (x,y) \in P \mid x \prec a \textup{ or } b \prec y \}
 \subseteq \s \subseteq
\{ (x,y) \in P \mid x \preceq a \textup{ or } b \preceq y \}.
\]
\end{conclusion}

\begin{proof}
The conclusion follows from the following implications
for a point $(x,y) \in P$:
\begin{center}
$ x \prec a \Longrightarrow (x,x) \in \s \Longrightarrow
 (x,y) \in \s $ \\
$b \prec y \Longrightarrow (y,y) \in \s \Longrightarrow
 (x,y) \in \s$ \\
$a \prec x \preceq y \prec b  \Longrightarrow
 (x,y) \in P \setminus \s$
\end{center}
\end{proof}

Given $a$ and $b$ in $X$, let 
\begin{align*}
H &= \{ (a,y) \in P \mid
  a \preceq y \preceq b \}, \text{ and} \\ 
V &= \{ (x,b) \in P \mid
  a \preceq x \preceq b \}.
\end{align*}

\begin{fact}\label{f:fact7}
If $a$ has no gap above, then $H \cap \s \subseteq 
 \{ (a,b) \}$. 
If $b$ has no gap below, then $V \cap \s \subseteq
 \{ (a,b) \}$.
\end{fact}

\begin{proof}
This follows immediately from the fact that $\s$ is
an open subset of $P$.
\end{proof}

\begin{fact}\label{f:fact8}
If $a$ has a gap above, then either $H \cap \s = H$
or $H \cap \s \subseteq \{ (a,b) \}$.
\end{fact}

\begin{proof}
Assume the contrary.  Then $(a,a) \notin \s$ and
there is $\beta$ such that $a \prec \beta \prec b$
and $(a, \beta) \in \s$.  Consider two cases.
First, assume that $\beta = \pred b$.  In this
case, since $b$ has a gap below, we also have
available the assumption that $(b,b) \notin \s$.
(See the comments after Fact~\ref{f:fact5}.)  Since
$(a,a) \notin \s$ and $(b,b) \notin \s$, both
$\pi(E_{[0, a]})$ and $\pi(E_{[b, 1]})$ are non-zero.
Furthermore, for any $g \in A$, $\supp E_{[0, a]}
g E_{[b, 1]} \subseteq \s$.  Indeed, if
$(x,y) \in \supp E_{[0, a]} g E_{[b, a]}$, then
$x \preceq a$ and $b \preceq y$.  If either
$x \prec a$ or $b \prec y$, the $(x,y) \in \s$.
If $x=a$ and $y=b$, then $(a,b) = (x,y) \in P$;
since $(a, \beta) \in \s$, we also have $(a,b) \in \s$.
Fact~\ref{f:fact4} implies that $\pi$ is not a nest
representation, a contradiction.

In the alternative case, $\beta \prec b$ and there is
$\barb$ such that $\beta \prec \barb \prec b$ and
$\barb$ has a gap below.  Since $(\barb, \barb) \notin \s$,
$\pi(E_{[\barb, 1]}) \ne 0$.  As before,
$\pi(E_{[0, a]}) \ne 0$.  Let $(x,y) \in
\supp E_{[0, a]} g E_{[\barb, 1]}$, where $g$ is any
element of $A$.  If $x \prec a$ then $(x,y) \in \s$.
If $x =a$ then $\beta \prec \barb \preceq y$, whence
$(x,y) = (a,y) \in \s$.  Once again, Fact~\ref{f:fact4}
yields a contradiction
\end{proof}

\begin{fact}\label{f:fact9}
If $b$ has a gap below, then either $V \cap \s =V$
or $V \cap \s \subseteq \{ (a,b) \}$.
\end{fact}

\begin{proof}
The idea behind the proof is essentially the same as in
the proof of Fact~\ref{f:fact8}.  This time, if the
conclusion does not hold, then
$(b,b) \notin \s$ and there is $\alpha$ 
such that $a \prec \alpha \prec b$ and $(\alpha, b) \in \s$.
If $\alpha = \suc a$ then take $\bara = a$; otherwise, 
take $\bara$ so that $a \prec \bara \prec \alpha$ and
$\bara$ has a gap above.  Now apply Fact~\ref{f:fact4}
to $E_{[0, \bara]}$ and $E_{[b, 1]}$.
\end{proof}

\begin{conclusion}\label{c:concB}
Conclusion~\ref{c:concA} and Facts~\ref{f:fact8}
and \ref{f:fact9} imply that $\s$ has one of
the following two forms:
\begin{align*}
\sab &= \{ (x,y) \in P \mid x \prec a \text{ or }
  b \prec y \} \quad \text{or} \\
\tab &= \sab \cup \{ (a,b) \}
\end{align*}
\end{conclusion}

Note that the latter is a possibility only if
$(a,b) \in P$ and $\tab$ is open.

\begin{fact}\label{f:fact10}
If $a$ has a gap above and $b$ has a gap below and
either $\s = \sab$ with $(a,b) \notin P$ or
$\s = \tab$, then $\pi$ is not a nest representation.
\end{fact}

\begin{proof}
Apply Fact~\ref{f:fact4} to $E_{[0,a]}$ and $E_{[b, 1]}$.
\end{proof}

This effectively ends the proof of Theorem~\ref{t:miker}.  
If $\pi$
is a nest representation, then $\s$ is one of the
ideal sets listed in Theorem~\ref{t:misets}.  Since
these are all meet irreducible, we have proven that
$\ker \pi$ is meet irreducible.
 \end{proof}

\section{Nest Representations with Atoms} \label{s:atomic}

In this section we give a condition on a nest representation which
guarantees that $\ker \pi$ is meet irreducible.  This condition
requires that $\pi$ be a \star representation on the diagonal of the
strongly maximal TAF algebra.  As pointed out in the introduction, any 
nest representation is similar to one with this property;
consequently, we assume throughout this section that the restriction
of $\pi$ to $D$ is a \star representation.

Recall that, in a von Neumann algebra 
$\DD$, a projection $E$ is said to be an
 \textit{atom} if $E$ majorizes no proper
(nonzero) subprojection.  
$\DD$ is \textit{atomic} if, for any projection $P \in \DD$,
$P = \vee \{ E \in \DD: E \text{ is an atom and } E \leq P \}$.

If $\pi$ is a nest representation
of a strongly maximal TAF algebra $A$ with diagonal $D$
such that the von Neumann algebra $\pi(D)'' $
contains an atom, then $\ker \pi$ is meet irreducible.
This is established in Theorem~\ref{t:atomic}, for which we give two
proofs.  The first proof depends on Theorem~2.1 in \cite{dhhls}; the
alternative proof is independent of this theorem.
 Both proofs require the fact (established in Proposition~\ref{p:zero1})
 that if $\pi(D)''$ contains an 
atom, then it is an atomic von Neumann algebra.
  The alternative proof can be read
immediately after Proposition~\ref{p:zero1}; from this point on it
uses the inductivity of ideals  rather than the spectrum 
characterization of meet irreducible ideals from \cite{dhhls}.
  In fact, a reader willing 
to assume that $\pi(D)''$ is  atomic can read the
alternative proof to Theorem~\ref{t:atomic} immediately after
Lemma~\ref{l:atom2}.  This provides a much shorter route to a
somewhat weaker theorem.  The alternative proof can be found at the
end of this section.

If $\pi(D)'' $ is atomic, then Corollary~\ref{c:dec} implies that
$\lat \pi(A)$ is a purely atomic nest.  The reverse implication is
false: Example I.3 in \cite{jojp95} provides an example of a
\star extendible
representation of a standard  limit algebra for which 
$\lat \pi(A) = \{0,I\}$ and $\pi(D) $ is weakly dense in a continuous
masa. Since $\ker \pi = \{0\}$ in this example, $\ker \pi$ is meet
irrreducible; we conjecture that $\ker \pi$ is meet irreducible any
time that $\lat \pi(A)$ contains an atom.  Proposition \ref{p:zero1}
establishes a dichotomy: either $\pi(D)''$ is  atomic  or else
$\pi(D)''$ contains no atoms at all.  This raises the natural
question: is the same dichotomy valid for $\lat \pi(A)$?
  This dichotomy certainly does not hold
for general nest representations (those which are not \star
representations on the diagonal); the failure is a consequence of the
similarity theory for nests.  The proof of Theorem 2.4 in
\cite{dhhls} applied to a refinement algebra provides an example of a
nest representation $\pi$ of a TAF agebra whose nest is purely
atomic, order isomorphic to the Cantor set, and such that the atoms are
ordered as the rationals.  There exist nests which are not purely
atomic but which are order isomorphic to the Cantor nest and which
have (rank
one) atoms  ordered as the rationals; the similarity theorem for
nests \cite{krd84} gives the existence of an invertible operator which
carries the invariant subspace nest for $\pi$ onto a nest of the
second type.  The composition of $\pi$ with this similarity yields a
nest representation of a strongly maximal TAF algebra (the refinement
algebra) which has atoms but is not purely atomic.

As before, we let $X$ denote the spectrum of the diagonal $D$ of a
TAF algebra $A$.  This is a zero dimensional topological space and
the clopen sets form a basis for the topology.  When $e$ is a 
projection in $D$, we let $\ehat$ denote the spectrum of $e$
in $X$ (i.e., the support set of $e$ viewed as an element of $C(X)$).

Now suppose that $\pi$ is a \star representation of the diagonal
$D$ of a TAF algebra and let $E$ be the spectral measure associated
with $\pi$.  $E$ is a regular, projection valued measure defined on
the Borel sets of $X$ which ``agrees'' with $\pi$ on clopen subsets in
the sense that $E(\ehat) = \pi(e)$, where $e$ is any projection in $D$
and $\ehat$ is its support in $X$.  If $\DD$ is the von Neumann
algebra generated by $\pi(D)$, then any projection $P$ in $\DD$ is of
the form $E(S)$, where $S$ is a Borel subset of $X$.

When $S$ is a singleton $\{x\}$ we shall write $E_x$ in place of
$E(\{x\})$.  If $e_n$ is a decreasing sequence of projections in $\DD$
such that $\cap \,\ehat_n = \{x\}$, then, by the regularity of $E$,
$\ds E_x = \wedge \pi(e_n)$.  In particular, $\wedge \pi(e_n) =
\wedge \pi(f_n)$ for any two decreasing sequences of projections
in $\DD$ with $\cap \,\ehat_n = \{x\}$ and $\cap \,\fhat_n = \{x\}$.

If there is a projection $e$ in $\ker \pi$ with $x \in \ehat$,
then clearly $E_x = 0$.  In this case, if $\ehat_n$ is any decreasing
sequence of clopen sets with $\cap \,\ehat_n = \{x\}$, then, since
$\{\ehat_n\}$ is a neighborhood basis for $x$, we have
$\pi(e_n) = 0$ for all large $n$.  The alternative is that
$\pi(e) \neq 0$ for any projection with $x \in \ehat$; in particular,
for any decreasing sequence $e_n$ with $\cap \,\ehat_n = \{x\}$,
$\pi(e_n) \neq 0$ for all $n$.  The projection $E_x = \wedge \pi(e_n)$
may or may not be $0$; it is, however, independent of the choice of
decreasing clopen sets with intersection $\{x\}$.

Note also that if $x,y \in X$ and $x \neq y$, then $E_xE_y = 0$.

\begin{lemma} \label{l:atom2} Let $D$ be the diagonal of a 
TAF algebra; let $\pi \colon D \to \BB(\HH)$ be a
\star representation; let $E$ be the spectral measure for
$\pi$; and let $\DD$ be the von Neumann algebra generated
by $\pi(D)$.  For any $x \in X$, if $E_x \neq 0$ then
$E_x$ is an atom of $\DD$.  Conversely, if $E_0$ is an atom for
$\DD$, then there is a unique element $x \in X$ such that $E_0 = E_x$.
\end{lemma} 

\begin{proof}
Any projection in $\DD$ has the form $E(S)$ for some Borel subset
$S$ of $X$.  Given $x \in X$, if $x \in S$ then $E_x \leq E(S)$ and
if $x \notin S$ then $E_x E(S) = 0$.  This shows that when 
$E_x \neq 0$, it is an atom of $\DD$.

Now suppose that $E_0$ is an atom of $\DD$.  Let $S$ be such that
$E_0 = E(S)$.  It is evident that there is at most one point
$x \in S$ such that $E_x \neq 0$; we need to prove the existence of
such a point.

Since $X$ is a Cantor set, we can find, for each $n \in \mathbb{N}$, 
$2^n$ disjoint clopen sets $\ehat^n_k$, $k = 1, \dots , 2^n$,
whose union is $X$, with the further property that any decreasing
sequence of these sets has one-point intersection.
Since $E_0$ is an atom of $\DD$, for each $n$ there is a unique
integer $k_n$ in $\{1, \dots, 2^n \}$ such that 
$E_0 \leq \pi(\ehat^n_{k_n})$.  Let $x$ be such that
$\ds \cap_n \,\ehat^n_{k_n} = \{x\}$.  Clearly, $E_0 \leq E_x$.
But $E_x$ is an atom when it is non-zero; hence $E_0 = E_x$.
The uniqueness of $x$ follows immediately form the orthogonality
of $E_x$ and $E_y$ when $x \neq y$.
\end{proof}

\begin{notation}
Let $\pi \colon A \to \BB(\HH)$ be a representation of a TAF algebra 
$A$ (which acts as a \star representation on the diagonal $D$).  For
$\xi \in \HH$ and $x \in X$ let $\MM_{\xi}$ denote the smallest
$\pi$-invariant subspace which contains $\xi$ and $\MM_x$ denote
the smallest $\pi$-invariant subspace which contains $\Ran E_x$. 
\end{notation}

Since the linear span of the $D$-normalizing partial isometries is 
dense in $A$, $\MM_{\xi}$ is the closed linear span of
$\{ \pi(v)\xi \mid v \in N_D(A) \}$ and $\MM_x$ is the closed linear
span of $\{ \pi(v)\xi \mid v \in N_D(A), \xi \in E_x \}$.

We now need to investigate the manner in which $D$-normalizing partial
isometries act on atoms of $\DD$.  Throughout the remainder of this
section,  $A$ will denote a TAF algebra (we will add the hypothesis
that $A$ is strongly maximal later); $\pi$ will denote a nest
representation of $A$ acting on a Hilbert space $\HH$; $D$ will denote
the diagonal of $A$ (with spectrum $X$); and $\DD = \pi (D)''$.

\begin{lemma} \label{l:npi} Let $v \in N_D(A)$ and let $x \in X$. 
 If $x \notin \wh{v^*v}$, then
$\pi (v)E_x = 0$.
 If $x \in \wh{v^*v}$, there exists $y \in X$ 
such that $(y, x) \in \vhat$ and  $\pi (v)E_x = E_y \pi(v)$.
In particular, if $y \neq x$, $\Ran \pi (v)E_x \perp \Ran E_x $.
\end{lemma}

\begin{proof}  Let $\{ e_n\}$ be a decreasing sequence 
of projections in $D$ with $\cap \,\ehat_n = \{ x\}$.
If $x \notin \wh{v^*v}$, then $ve_n = 0$ for large $n$, 
in which case
$\pi(v) E_x = \pi (v) \pi(e_n)E_x = \pi(ve_n)E_x = 0$.

 Now suppose that $x \in \wh{v^*v}$ and let $y$ 
be such that $(y, x) \in \vhat$. 
With $e_n$ as above, let $f_n = ve_nv^*$, so that 
$\cap \,\fhat_n = \{ y\}$ and $\wedge \pi (f_n) = E_y$.
Then $ ve_n = f_nv$ for  large $n$; hence
$\pi (v) \pi (e_n) = \pi (f_n) \pi (v) $ and, taking
strong limits,
$\pi(v) E_x =  E_y \pi (v)$.
It follows that $\Ran (\pi (v) E_x) \subset \Ran E_y$ 
and, hence, that $E_x$ and $\pi (v)E_x$
have orthogonal ranges when
$y \neq x$.

\begin{remark}
If $v \in N_D(A)$ and $(y,x) \in \vhat$, then 
$\pi(v)E_x = E_y \pi(v) E_x$.  It is possible that
$\pi(v)E_x = 0$ even when $\pi(v) \neq 0$ and 
$E_x \neq 0$.
\end{remark}

\end{proof}
\begin{lemma} \label{l:rkone} 
Let $x \in X$.
If $E_x \neq 0$, then $E_x$ is a rank-one atom.
\end{lemma}

\begin{proof} 
Let $x \in X$ and assume that $\Ran E_x$ contains two
unit vectors $\xi$ and $\zeta$ such that $\xi \perp \zeta$.
Let $v \in N_D(A)$.  Since $\pi(v)E_x = \pi(v)\pi(e)E_x$ for
any projection $e$ for which $x \in \ehat$, we may, by a 
suitable restriction, reduce to considering two cases:
when $\vhat$ is contained in the diagonal $\{(z,z) \mid
z \in X \}$ of $G$ and when $\vhat$ is disjoint from the
diagonal.  When $\vhat$ is contained in the diagonal,
$v$ is a projection in $D$ and $\pi(v)$ either dominates $E_x$
or is orthogonal to $E_x$.  In this case, either 
$\pi(v)\xi = \xi$ and $\pi(v)\zeta = \zeta$ or $\pi(v)\xi = 0$ 
and $\pi(v)\zeta = 0$.  When $\vhat$ is disjoint form the diagonal
of $G$, there is $y \in X$ such that $y \neq x$ and
$(y,x) \in \vhat$.  In this case, $\pi(v)\xi \in \Ran E_y$
and $\pi(v) \zeta \in \Ran E_y$.  In particular,
$\pi(v)\xi \perp \zeta$ and $\pi(v) \zeta \perp \xi$ for all
$v \in N_D(A)$.

It now follows that $\xi \in \MM_{\xi}$ and 
$\xi \notin \MM_{\zeta}$ while $\zeta \in \MM_{\zeta}$
and $\zeta \notin \MM_{\xi}$.  But $\pi$ is a nest representation
and $\MM_{\xi}$ and $\MM_{\zeta}$ are $\pi$-invariant; hence one must
contain the other.  Thus, the rank of $E_x$ is at most $1$.
\end{proof}

\begin{corollary} \label{c:2pis}  Let 
$u,w \in N_D(A)$ and assume that $u$ and $w$ 
have a common subordinate.
Let $(y, x)$ be a point in $P$ such that  $(y, x) \in \uhat \cap
\what$.  Then $\pi(u)E_x = \pi(w)E_x$ and,
if $E_x$, $E_y \neq 0$, $ \Ran \pi(u)E_x = \Ran E_y$.
\end{corollary}

\begin{proof} Let $e$ be a (nonzero) projection in 
$D$ such that $e \leq \min \{ u^*u, w^*w\}$ 
and $x \in \ehat$. (One could just take
$e = u^*uw^*w$.)  Then $ue = we$.  Since $\pi(e) $ dominates $E_x$,
\[ 
\pi(u)E_x = \pi(u)\pi(e)E_x = \pi(ue)E_x = \pi(we)E_x 
= \pi(w)\pi(e)E_x = \pi(w)E_x.
\]

For the second assertion, suppose $E_x$ and $E_y \neq 0$.  
Now by Lemma~\ref{l:npi}, for any 
 $v \in N_D(A)$, either $\pi(v)E_y$ is zero or else 
$\pi(v)E_y$ is contained in $\Ran E_z$ for some $z \in X$
with $(z, y) \in \vhat$.  But as $y \prec x$, 
$\Ran \pi(v)E_y \perp \Ran E_x$. 
We have  seen that $\MM_y$, the smallest 
$\pi$-invariant subspace containing
$\Ran E_y$, is disjoint from
$\Ran E_x$; since $\pi$ is a nest representation, 
it follows that $\MM_y \subset \MM_x$.  Thus, for some 
$v$, $\Ran \pi(v)E_x \cap \Ran E_y \neq (0)$. 
 Since the ranges of $E_x$ and $E_y$
are one-dimensional, $\Ran \pi(v)E_x = \Ran E_y$.  
For such a $v$ we have $(y, x) \in \vhat$.
By the first paragraph, $v$ can be replaced by any 
$u$ with $(y, x) \in \uhat$.
\end{proof}

\begin{proposition} \label{p:zero1}
 $E_0 = \vee \{ E_x \mid x \in X \}$
is either $0$ or $I$.
\end{proposition}

\begin{proof} 
Let $E_1 = I - E_0$.  We shall show that both $E_0\HH$
and $E_1\HH$ are invariant under $\pi$.  Since $\pi$
is a nest representation, this means that one of
them must be zero.

 If $E_0\HH$ is not invariant, then for some $ v \in N_D(A)$
and $\xi \in E_0\HH$, 
$E_1\pi(v)\xi \neq 0$. 
Now $\xi = \sum \{ E_x \xi \mid x \in X \}$, so there exists 
$x \in X$ with $E_1\pi(v)E_x\xi \neq 0$. 
By Lemma~\ref{l:npi},
$\pi(v)E_x = E_y\pi(v)E_x$, for some $y$.  Since $\pi(v)E_x \neq 0$ 
it follows that $E_y \neq  0$; hence $E_y$ is an atom
in $\DD$.  As $E_1$ majorizes no atoms, $E_1E_y = 0$.  
On the other hand, $0 \neq E_1\pi(v)E_x
=E_1E_y\pi(v)E_x$, so that $E_1E_y \neq 0$.  
This contradiction shows that $E_0\HH$ is invariant.

If $E_1 \HH$ is not invariant, there exists 
 a vector $\xi \in E_1\HH$ and a $D$-normalizing partial
isometry $v$ in $A$ such that $E_0\pi(v)E_1\xi \neq 0$. 
Since $E_0$ is the sum of the atoms it majorizes, 
$E_y\pi(v)E_1\xi \neq 0$  
for some $y \in X$. 
By Lemma~\ref{l:npi}, there is an element $x\in X$ with
$E_y\pi(v) = E_y\pi(v)E_x$,  whence $E_y\pi(v)E_xE_1 \neq 0$.  
In particular, $ E_xE_1 \neq 0$, which
contradicts the fact that $E_1$ majorizes no atoms.  
Thus, $E_1\HH$ is also invariant. 
\end{proof}

\begin{remark} Let $\DD$ be the von Neumann algebra generated 
by $\pi(D)$, where $D$ is the diagonal of the
TAF algebra $A$.  According to Proposition~\ref{p:zero1}, 
if $\pi \colon A \to \BB(\HH)$ is a nest representation, then
either $\DD$ is generated by its atoms (that is, $\DD$ is 
purely atomic), or else it has no atoms.
(This, of course, presumes that the restriction of $\pi$ to
$D$ is a \star representation.)
\end{remark}

\begin{lemma} \label{l:orbit} 
If $x$ and $y$ are two points of
$X$ such that $E_x$ and $E_y$ are both nonzero, then 
$x$ and $y$ belong to the same orbit in $X$.
\end{lemma}

\begin{proof} 
If $x$ and $y$ are not in the same orbit, it follows
 that $\Ran E_x \perp \MM_y$
and $\Ran E_y \perp \MM_x$  (Lemma~\ref{l:npi}).
Since $\Ran E_y \subseteq \MM_y$ and
$\Ran E_x \subseteq \MM_x$, 
$\MM_x$ and $\MM_y$ are not lineraly ordered, a contradiction.
\end{proof}

\begin{lemma} \label{l:inter}
Assume, further, that
 $A$ is a strongly maximal TAF algebra. 
If, for some $x \in X$, $E_x \neq 0$, then 
$J = \{ z \mid E_z \neq 0 \}$ is an interval in the orbit of $x$.
\end{lemma}

\begin{proof}
If $E_x \neq 0$, Lemma~\ref{l:orbit} implies that $J$ is 
contained in the orbit of $x$.  If $E_z = 0$
for all $z \neq x$, we are done.  Suppose then, that
 $E_y \neq 0$  
for some $y \neq x$ in the orbit of $x$.  Without loss of generality
we may assume that  $y \prec x$.  Let $v \in N_D(A)$ be
such that $(y,x) \in \vhat$.  Then
 $\MM_y \subset \MM_x$ and    
$\Ran \pi(v)E_x = \Ran E_y$
(Corollary~\ref{c:2pis}).
In particular, $\pi(v)E_x \neq 0$.

Let $z$ be a point in the orbit of $x$ with 
$y \prec z \prec x$.  Since $A$ is strongly maximal, there exist
$u, w \in N_D(A)$ with 
$ (y, z) \in \uhat$ and $(z, x) \in \what$.
Now $\pi(v)E_x = \pi(uw)E_x 
= \pi(u)\pi(w)E_x$; hence $\pi(w)E_x \neq 0$.
As $\Ran \pi(w)E_x \subset \Ran E_z$, 
it follows that $E_z \neq 0$ and, hence, $z \in J$.  
This shows that $J$ is an interval.
\end{proof}

\begin{corollary} \label{c:masa} 
If $\DD$ has an atom
and $J$ is the interval obtained in Lemma~\ref{l:inter}, then
$\vee \{E_x \mid x\in J\} = I$ and
 $\DD$ is a masa in $\BB(\HH)$.

\end{corollary}
\begin{proof} The first assertion follows from 
Lemma~\ref{l:atom2} and Proposition~\ref{p:zero1}, 
since all atoms have the form $E_x$, for some
$x \in X$.
 Since the set $\{ E_x \mid x\in X \}$ is a collection of 
commuting, rank-one atoms whose ranges span $\HH$, the
von Neumann algebra which they generate is a masa in $\BB(\HH)$.
\end{proof}

\begin{theorem} \label{t:atomic} Let $A$ be a strongly maximal 
TAF algebra and $\pi \colon A \to \BB(\HH)$ be a
nest representation for which the von Neumann 
algebra $\DD$ generated by $\pi(D)$ contains an atom.
Then the kernel of $\pi$ is a meet irreducible ideal in $A$.
\end{theorem}

\begin{proof}
By hypothesis, $\DD$ contains an atom, necessarily of the form
$E_x$, for some $x \in X$.
By Lemma~\ref{l:inter}, there is a nonempty interval $J$ 
in an orbit in $X$ with the following property:
 for any $D$-normalizing partial isometry $v$,
 $\pi(v) \neq 0$ if, and only, if  
   $J \times J$ intersects $\vhat$. 
In other words, the complement of the spectrum of the ideal
$\ker \pi$ contains $(J \times J)\cap P$.
But the complement is a closed set, so it contains 
$\overline{(J \times J)\cap P}$.  On the other hand,
if $v$ is a $D$-normalizing partial isometry with $\vhat$ 
disjoint from $\overline{(J \times J) \cap P}$,
then $\pi(v) = 0$.  Thus,
\[ 
\wh{\ker \pi} = P \setminus \overline{(J \times J) \cap P}.
\]
 By  \cite[Theorem 2.1]{dhhls},  $\ker \pi$ is a meet-irreduclble ideal.
\end{proof}

\begin{corollary} \label{c:dec} Let $A$, $\pi$, and $J$ 
be as above and let $\NN = \lat \pi (A)$.
\begin{theoremlist}
\item If $F$ is decreasing subset of $J$, then 
$P = \sum \{ E_x \mid x \in F\} \in  \NN$.
On the other hand,
 if $P\in  \NN $, then $F = \{ x \mid 0 \neq E_x \leq P \}$ is a decreasing
subset of $J$.  This correspondence between decreasing subsets of $J$
and projections in $\NN$ is a bijection.
\item $\DD = \alg \NN \cap (\alg \NN)^*$ 
equals the von Neumann algebra generated by $\NN$.
\end{theoremlist}

\end{corollary}
\begin{proof} The first assertion is clear.  For the second,
first note that, since $\DD$ is a masa, 
$\DD = \alg \NN \cap (\alg \NN)^*$. 
Let $x\in J$, let $P = \sum \{E_y \mid y \in J \text{ and } 
y \preceq x\}$, and let
$P_{-} = \sum \{E_y \mid y \in I \text{ and } y \prec x\}$. 
Then $P_{-}$ is the immediate predecessor of $P$ in $\NN$ and
$E_x = P - P_{-}$. 
Thus, every atom from $\DD$, and hence $\DD$ itself, is
contained in the von Neumann algebra generated by
$\NN$.  The reverse inclusion is obvious.
\end{proof}

There is an alternate proof for Theorem~\ref{t:atomic} based on
a presentation for $A$ rather than on the spectrum of $A$ and
Theorem 2.1 in \cite{dhhls}.  This proof is dependent only on the
preliminary results through Proposition~\ref{p:zero1}; in fact, if the
reader is willing to assume that $\pi(D)'' $ is purely
atomic, then only Lemma~\ref{l:atom2} is needed.  In the alternate
proof, we view $A$ as the union of an ascending chain of 
 subalgebras each of which is star extendibly isomorphic to a maximal
triangular subalgebra of a finite dimensional \cstar algebra.
Also we may assume that a system of matrix
units for each $A_k$ has been selected in such a way that each matrix
unit in $A_k$ is a sum of matrix units in $A_{k+1}$; this gives a
matrix unit system for $A$.  Another fact from the lore of direct
limit algebras that we need is that ideals are inductive: if $I$ is
an ideal in $A$, then $I$ is the closed union of the
ideals $I_k = I \cap A_k$ in $A_k$.

\begin{proof}[Alternate Proof of Theorem~\ref{t:atomic}]
Assume that $\pi$ is a nest representation but that $\ker \pi$
is not meet irreducible.  Let $I$ and $J$ be two ideals in
$A$ such that $I \cap J = \ker \pi$ and $I \cap J$ differs from
both $I$ and $J$.  By the inductivity of ideals,
there exist matrix units $u_I \in I \setminus J$ and
$u_J \in J \setminus I$.  These matrix units must lie in some
$A_k$; since we may replace the sequence $A_k$ by a subsequence, 
 we may assume that $u_I$ and $u_J$ lie in $A_1$.  

Since $u_I \notin J$, $\pi(u_I) \neq 0$. By Proposition~\ref{p:zero1},
\[
\pi(u_I) = \sum_{x,y} E_x \pi(u_I) E_y,
\]
where the sum is taken over all pairs of atoms and is convergent
in the strong operator topology.  Consequently, there exist points
$x_I$ and $y_I$ in $X$ such that 
$E_{x_I} \pi(u_I) E_{y_I} \neq 0$.
For each $k$, let $e^I_k$ and $f^I_k$ be the unique diagonal
matrix units in $A_k$ such that $x_I \in \ehat^I_k$ and
$y_I \in \fhat^I_k$.  Since $u_I \in I_k$, 
$e^I_k u_I f^I_k \in I_k$, for all $k$.  On the other hand,
$\pi(e^I_k u_I f^I_k) \neq 0$, so $e^I_k u_I f^I_k \notin J_k$,
for all $k$.

        Let $x_J$ and $y_J$ in $X$ and 
$e^J_k$ and $f^J_k$ in $A_k$ be analogously
defined for the ideal $J$.  This time, $e^J_k u_J f^J_k \in J_k$ and
$e^J_k u_J f^J_k \notin I_k$.

        We shall show (after possibly reversing the 
roles if $I$ and $J$) that for infinitely many $k$,
$e^I_k A_k f^J_k \subseteq \ker \pi$ and 
$f^J_k A_k e^I_k \subseteq \ker \pi$.  This leads to 
a contradiction with the hypothesis that $\pi$ is
a  nest representation and so shows that $\ker \pi$ is
necessarily meet irreducible.  Indeed, from
\begin{align*}
\pi(e^I_k) \pi(A_k) \pi(f^J_k) &= 0, \quad \text{all $k$}, \\
\pi(f^J_k) \pi(A_k) \pi(e^I_k) &= 0, \quad \text{all $k$},
\end{align*}
it follows that
\begin{align*}
E_{x_I} \pi(A) E_{y_J} &=0, \\
E_{y_J} \pi(A) E_{x_I} &=0,
\end{align*}
and hence that $E_{x_I} \perp \MM_{y_J}$ and
$E_{y_J} \perp \MM_{x_I}$, where $\MM_{y_J}$ and
$\MM_{x_I}$ are the smallest $\pi$-invariant subspaces
containing $E_{y_J}$ and $E_{x_I}$ respectively.
But then $\MM_{y_J}$ and $\MM_{x_I}$ are not
related by inclusion, a contradiction.

        Each finite dimensional algebra $A_k$ is a direct sum
of $T_n$'s  and the matrix units $e^I_k$ and $f^I_k$ are 
in the same summand, as are $e^J_k$ and $f^J_k$.  If these two
summands differ, then $e^I_k A_k f^J_k = 0$ and
$f^J_k A_k e^I_k = 0$.  Should this occur for infinitely many
$k$, then we are done.  So we need consider only the case in
 which, for all $k$,
 all of $e^I_k$, $e^J_k$, $f^I_k$ and $f^J_k$ are in the same
summand in $A_k$.

If $e$ and $f$ are diagonal matrix units (minimal diagonal
projections) in  $A_k$, let $\m (e,f)$ be the matrix unit
in $C^*(A_k)$ with initial projection $f$ and final 
projection $e$ (if there is such a matrix unit).  If
$\m(e,f) \in A_k$, then $e \preceq f$ in the diagonal order
on minimal diagonal projections.  We shall need the following
property of ideals in $A_k$: if $e_1 \preceq e_2 \preceq
f_2 \preceq f_1$ and if $\m(e_2, f_2)$ is in an ideal, 
then $\m(e_1, f_1)$ is also in the ideal.

Since $e^I_k$ and $e^J_k$ are in the same $T_n$-summand of
$A_k$, they are related in the diagonal order.  By interchanging 
$I$ and $J$ and passing to a subsequence, if necessary, we may
assume that $e^I_k \preceq e^J_k$, for all $k$.  The facts
concerning the membership of $e^I_k u_I f^I_k$ and
$e^J_k u_J f^J_k$ in $I_k$ and $J_k$ may be rephrased as
\begin{align*}
\m(e^I_k, f^I_k) \in I_k &\text{ and } 
\m(e^I_k, f^I_k) \notin J_k, \\
\m(e^J_k, f^J_k) \notin I_k &\text{ and } 
\m(e^J_k, f^J_k) \in J_k. 
\end{align*}
As a consequence $f^I_k \prec f^J_k$, for all k.
(If $f^J_k \preceq f^I_k$ for some $k$, then
$e^I_k \preceq e^J_k \preceq f^J_k \preceq f^I_k$.
Since $\m(e^J_k, f^J_k) \in J_k$, we have
$\m(e^I_k, f^I_k) \in J_k$, a contradiction.)
But now,
\begin{align*}
\m(e^I_k, f^I_k) \in I_k \text{ and } f^I_k \prec f^J_k
&\Longrightarrow \m(e^I_k, f^J_k) \in I_k \\
\m(e^J_k, f^J_k) \in J_k \text{ and } e^I_k \preceq e^J_k
&\Longrightarrow \m(e^I_k, f^J_k) \in J_k.
\end{align*}
Thus $\m(e^I_k, f^J_k) \in I_k \cap J_k \subseteq \ker \pi$;
hence $e^I_k A_k f^J_k \subseteq \ker \pi$.
Also, since $e^I_k \preceq f^I_k \prec f^J_k$,
$f^J_k  A_k e^I_k = \{0\} \subseteq \ker \pi$.
As pointed out earlier, this implies that $\pi$ is not
 a nest representation; so, when $\pi$ is a nest 
representation with an atomic lattice, $\ker \pi$ is
meet irreducible.
\end{proof}

\providecommand{\bysame}{\leavevmode\hbox to3em{\hrulefill}\thinspace}

\end{document}